\documentclass[11pt]{article}
\usepackage{amsmath}
\usepackage{amssymb}
\usepackage{latexsym}
\usepackage{graphicx}
\newcommand{\ds}{\displaystyle}
\newtheorem{theorem}{\bf Theorem}
\newtheorem{lemma}{\bf Lemma}
\newtheorem{coro}{\bf Corollary}
\newtheorem{remark}{\bf Remark}

\begin{document}

\baselineskip 6mm

\vspace*{10mm}

\begin{center}
{\LARGE\bf Application of Tauberian Theorem to the Exponential\\[2mm]
Decay of the Tail Probability of a Random Variable}\\

\vspace*{10mm}

Kenji Nakagawa \\
Department of Electrical Engineering, \\
Nagaoka University of Technology, \\
Nagaoka, Niigata 940-2188, Japan \\
E-mail nakagawa@nagaokaut.ac.jp \\
\end{center}

\vspace*{5mm}

\begin{abstract}
We give a sufficient condition for the exponential decay 
of the tail probability of a non-negative random variable. 
We consider the Laplace-Stieltjes transform of the 
probability distribution function of the random variable. 
We present a theorem, according to which if the abscissa 
of convergence of the LS transform is negative finite and 
the real point on the axis of convergence is a pole of the 
LS transform, then the tail probability decays exponentially. 
For the proof of the theorem, we extend and apply 
so-called a finite form of Ikehara's complex Tauberian 
theorem by Graham-Vaaler. 

\end{abstract}


\vspace*{5mm}

\section{Introduction}
The purpose of this paper is to give a sufficient condition 
for the exponential decay of the tail probability of a 
non-negative random variable. For a non-negative random 
variable $X$, $P(X>x)$ is called the {\it tail probability} 
of $X$. The tail probability {\it decays exponentially} if 
the limit 
\begin{eqnarray}
\lim_{x\to\infty}\frac{1}{x}\log P(X>x)
\end{eqnarray}
exists and is a negative finite value. 

For the random variable $X$, the probability distribution 
function of $X$ is denoted by $F(x)= P(X\le x)$ and 
the Laplace-Stieltjes transform of $F(x)$ is denoted 
by $\varphi(s)=\int_0^\infty e^{-sx}dF(x)$. We will give a sufficient condition 
for the exponential decay of the tail probability $P(X>x)$ 
based on analytic properties of $\varphi(s)$. 

In \cite{nak2}, we obtained a result that 
the exponential decay of the tail probability $P(X>x)$ 
is determined by the singularities of $\varphi(s)$ on 
its axis of convergence. In this paper, we investigate 
the case where $\varphi(s)$ has a pole at the real point 
of the axis of convergence, and reveal the relation between 
analytic properties of $\varphi(s)$ and the exponential 
decay of $P(X>x)$. 

The results obtained in this paper will be applied to 
queueing analysis. In general, there are two main 
performance measures of queueing analysis, one is the 
number of customers $Q$ in the system and the other is 
the sojourn time $W$ in the system. $Q$ is a discrete 
random variable and $W$ is a continuous one. It is 
important to evaluate the tail probabilities $P(Q>q)$ 
and $P(W>w)$ for designing the buffer size or link capacity 
in communication networks. Even in the case that 
the probability distribution functions $P(Q\le q)$ or 
$P(W\le w)$ cannot be calculated explicitly, their 
generating functions $Q(z)=\sum_{q=0}^\infty P(Q=q)z^q$ 
or $W(s)=\int_0^\infty e^{-sw}dP(W\le w)$ can be obtained 
explicitly in many queues. Particularly, in M/G/1 queue, 
$Q(z)$ and $W(s)$ are given explicitly by Pollaczek-Khinchin 
formula \cite{kle}. 
So, in this paper, we assume that we have the explicit 
form of a generating function and then investigate 
the exponential decay of the tail probability based 
on the analytic properties of the generating function. 

Such kind of researches have been studied motivated 
by a requirement for evaluating a packet loss probability of 
a light tailed traffic in the packet switched network. 

An approach by the complex analysis is seen in \cite{fal}. 
A sufficient condition is given in \cite{fal} for the decay 
of the stationary probability of an M/G/1 type Markov 
chain with boundary modification and the result is applied 
to MAP/G/1 queue. Let $\pi=(\pi_n)$ denote the stationary 
probability of an M/G/1 type Markov chain with boundary 
modification, and $\pi(z)=\sum_n\pi_nz^n$ the probability 
generating function of $\pi$ with the radius of convergence 
$r>0$. In \cite{fal}, they proved a theorem that if $z=r$ 
is a pole of order $1$ and is the only singularity on 
the circle of convergence $|z|=r$, then there exists 
$K>0$, $\tilde{r}>r$ such that $\pi_n=Kr^{-n}+O(\tilde{r}^{-n})$.
Our Theorem 1 below is an extension of this theorem in \cite{fal}. 

In \cite{duf}, for the stationary queue length $Q$ of 
a queueing system which satisfies some large deviations 
conditions, it is shown that the $P(Q>n)$ decays as 
$P(Q>n)\simeq\psi\exp(-\theta n)$ with a positive constants 
$\psi$ and $\theta$. In \cite{gly}, it is shown that the stationary waiting 
time $W$ and queue length $Q$ decay exponentially for a broad 
class of queues with stationary input and service. 
In \cite{neu1},\cite{neu2}, the 
stationary distribution of M/G/1 or G/M/1 type Markov chains are 
deeply studied. In \cite{tak}, the tail of the waiting time in PH/PH/c 
queue is investigated. In \cite{liq}, a sufficient condition 
is given for the stationary probability of a Markov chain of 
GI/G/1 type to be light tailed. 

This paper is organized as follows. In section II, an 
application of our results to queueing analysis is 
presented. In section III, an example of a random variable 
is given whose tail probability does not decay exponentially. 
In section IV, some Tauberian theorems are introduced 
and the relation to our problems is stated. In section V, 
some lemmas are given and our main theorem is proved in 
the case of a pole of order 2. In section VI, the statement 
of lemmas and theorems are presented for a pole of arbitrary 
order. Finally, we summarize out results in section VII. 

Throughout this paper, we use the following symbols. 
${\mathbb C}$, ${\mathbb R}$, ${\mathbb Z}$, $\Re$ denote 
the set of complex numbers, real numbers, integers, and the 
real part of a complex number, respectively.

\section{Application to Queueing Analysis}
The author already had results on the exponential decay 
of the tail probability for a discrete random variable 
\cite{nak1}, and those for a general random variable \cite{nak2}.
The main theorem in \cite{nak1} is as follows. 

\medskip

\begin{theorem}\label{theo:1}$\mbox{\rm\cite{nak1}}$ Let $X$ be 
a random variable taking non-negative integral values, and 
$f(z)$ be the probability generating function of $X$. The 
radius of convergence of $f(z)$ is denoted by $r$ and 
$1<r<\infty$ is assumed. If the singularities of $f(z)$ 
on the circle of convergence $|z|=r$ are only a finite 
number of poles, then 
\begin{eqnarray}
\lim_{n\to\infty}\frac{1}{n}\log P(X>n)=-\log r.
\end{eqnarray}
\end{theorem}

\bigskip

We can apply Theorem \ref{theo:1} to queueing analysis as 
follows. 

Consider, for example, the number of customers $Q$ in the steady 
state of M/D/1 queue with traffic intensity $\rho$. The 
probability generating function $Q(z)=\sum_{q=0}^\infty P(Q=q)z^q$ 
of $Q$ is given by Pollaczek-Khinchin formula \cite{kle}:
\begin{eqnarray}
Q(z)=\frac{(1-\rho)(z-1)\exp(\rho(z-1))}{z-\exp(\rho(z-1))}.
\end{eqnarray}
The radius of convergence of $Q(z)$ is equal to the 
unique solution $z=r>1$ of the equation 
$z-\exp(\rho(z-1))=0$. 
Since $Q(z)$ is meromorphic in the whole finite complex plane 
$|z|<\infty$, in particular, the singularities of $Q(z)$ 
on the circle of convergence $|z|=r$ are only a finite number 
of poles. Therefore, by Theorem \ref{theo:1}, we know that 
the tail probability $P(Q>q)$ decays exponentially as 
$q\to\infty$. 

Next, in \cite{nak2}, the exponential decay of the tail 
probability $P(X>x)$ is investigated for a general 
non-negative real valued random variable $X$. The main 
theorem in \cite{nak2} is as follows.

\medskip

\begin{theorem}\label{theo:2}$\mbox{\rm\cite{nak2}}$ Let 
$X$ be a non-negative random variable, and $\varphi(s)$ 
be the Laplace-Stieltjes transform of the probability 
distribution function of $X$. The abscissa of convergence 
of $\varphi(s)$ is denoted by $\sigma_0$ and $-\infty<\sigma_0<0$ 
is assumed. If the singularities of $\varphi(s)$ on the 
axis of convergence $\Re s=\sigma_0$ are only a finite number 
of poles, then we have
\begin{eqnarray}
\lim_{x\to\infty}\frac{1}{x}\log P(X>x)=\sigma_0.
\end{eqnarray}
\end{theorem}

\bigskip

Let us apply Theorem \ref{theo:2} to queueing analysis. 
In this case, however, the situation is somewhat different from 
that in the case of discrete random variable. 

Consider the sojourn time $W$ in M/D/1 queue with traffic 
intensity $\rho$. Writing $W(s)$ as the Laplace-Stieltjes 
transform of the probability distribution function 
$P(W\le w)$ of $W$, we have \cite{kle}
\begin{eqnarray}
W(s)=\frac{(1-\rho)s\exp(-s)}{s-\rho+\rho\exp(-s)}.\label{eqn:LSofW}
\end{eqnarray}
The abscissa of convergence of $W(s)$ is the unique negative 
solution $s=\sigma_0$ of the equation 
$s-\rho+\rho\exp(-s)=0$. 
We can see that the singularity of $W(s)$ on the axis of 
convergence $\Re s=\sigma_0$ is only a simple pole $s=\sigma_0$ 
\cite{nak3}. In fact, the location of the poles of $W(s)$ are 
shown in Figure 1.
\begin{figure}[t]
\begin{center}
\includegraphics[height=10cm]{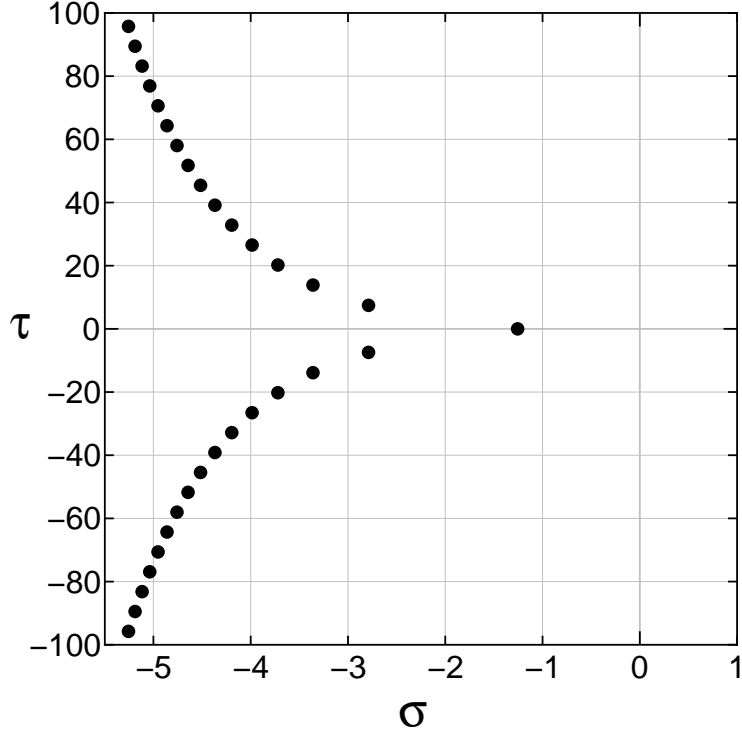}
\vspace*{-3mm}
\caption{The poles of $W(s)$ in (\ref{eqn:LSofW})\ 
with $\rho=0.5,\ s=\sigma+i\tau$}\label{fig:poles_LS}
\end{center}
\end{figure}
The abscissa of convergence is $\sigma_0=-1.26$ for 
$\rho=0.5$. Though, 
it is not easy to prove that $s=\sigma_0$ is the only 
singularity of $W(s)$ on the axis of convergence. In 
order to prove it, we need some theorems such as Rouch\'{e}'s 
theorem.

The statement of Theorem \ref{theo:1} and \ref{theo:2} are 
formally quite the same, but to verify that the assumption 
in Theorem \ref{theo:2} holds is more difficult than 
Theorem \ref{theo:1}. This is because of the difference of 
the convergence regions. In Theorem \ref{theo:1}, the 
boundary of the convergence region of a power series is a 
circle, which is a compact set, so if the probability generating 
function is meromorphic, then the singularities on the 
circle of convergence are necessarily a finite number of 
poles. On the other hand, in Theorem \ref{theo:2}, the axis 
of convergence is not a compact set, so we need some 
verification to see that the singularities on the axis of 
convergence are a finite number of poles. We want some simple 
sufficient condition to guarantee the exponential decay 
of the tail probability. We see, in fact, that the conclusion 
of Theorem \ref{theo:2} is really stronger than desired, 
so it may be possible to relax the assumption that the 
number of poles on the axis of convergence is finite. 

We have the following theorem. This is the main theorem in 
this paper.

\medskip

\begin{theorem}\label{theo:3} Let $X$ be a non-negative 
random variable, and $F(x)=P(X\le x)$ be the probability 
distribution function of $X$. Let 
\begin{eqnarray}
\varphi(s)=\int_0^\infty e^{-sx}dF(x),\ s=\sigma+i\tau\in{\mathbb C}
\end{eqnarray}
be the Laplace-Stieltjes transform of $F(x)$ and $\sigma_0$ 
be the abscissa of convergence of $\varphi(s)$. 
We assume $-\infty<\sigma_0<0$. If $s=\sigma_0$ is a pole 
of $\varphi(s)$, then we have
\begin{eqnarray}
\lim_{x\to\infty}\frac{1}{x}\log P(X>x)=\sigma_0.
\end{eqnarray}
\end{theorem}

\bigskip

\begin{remark} Since $F(x)$ is non-decreasing, $s=\sigma_0$ is 
a singularity of $\varphi(s)$ by Widder $\mbox{\rm\cite{wid}}$, 
p.$58$, Theorem $5$b. We assume in Theorem 
$\mbox{\rm\ref{theo:3}}$ that this singularity is a pole.
\end{remark}

\section{Example of Random Variable whose Tail Probability 
does not Decay Exponentially}

We show an example of a non-negative random variable whose 
tail probability does not decay exponentially, i.e., 
$x^{-1}\log P(X>x)$ does not have a limit \cite{nak2}. 

For any positive integer $h$, define a sequence 
$\{c_n\}_{n=0}^\infty$ by 
\begin{eqnarray}
\left\{
\begin{array}{ll}
c_0=0, & \\ \label{eqn:21}
c_n=c_{n-1}+h^{c_{n-1}}, & n=1,2,\ldots. 
\end{array}
\right.
\end{eqnarray}
We define a function $\gamma(x)$ by 
\begin{eqnarray}
\gamma(x)=h^{-c_{n}},\ \mbox{\rm for}\ c_n\le x<c_{n+1},\ n=0,1,\cdots.
\end{eqnarray}
For arbitrary $\sigma_0<0$, put $F^\ast(x)=1-e^{\sigma_0x}\gamma(x),
\ x\ge0$. 
We see that $F^\ast(x)$ is right continuous and non-decreasing 
with $F^\ast(0)=0$, $F^\ast(\infty)=1$, hence $F^\ast(x)$ is a 
distribution function. Let us define $X^\ast$ as a random 
variable with probability distribution function $F^\ast(x)$. 
We write $\varphi^\ast(s)$ as the Laplace-Stieltjes transform 
of $F^\ast(x)$. The following theorem shows that $X^\ast$ is 
an example of a random variable whose tail probability 
does not decay exponentially.

\medskip

\begin{theorem}\label{theo:4} $\mbox{\rm (see\,\cite{nak2})}$ Let $X^\ast$, 
$F^\ast(x)$ and $\varphi^\ast(s)$ be defined as above. 
Then, the abscissa of convergence of $\varphi^\ast(s)$ is $\sigma_0$, 
and we have 
\begin{align*}
\liminf_{x\to\infty}\frac{1}{x}\log P(X^\ast>x)
&\le\sigma_0-\log h\\
&<\sigma_0\\
&=\limsup_{x\to\infty}\frac{1}{x}\log P(X^\ast>x).
\end{align*}
All the points on the axis of convergence $\Re s=\sigma_0$ 
are singularities of $\varphi(s)$.
\end{theorem}

\bigskip

\begin{remark}
From Widder $\mbox{\rm\cite{wid}}$, p.$44$, Theorem $2.4$e, 
$\limsup_{x\to\infty}x^{-1}\log P(X>x)=\sigma_0$ holds 
under the condition $\sigma_0<0$ and no other condition is 
necessary. Meanwhile, Theorem $\ref{theo:4}$ implies 
that some additional condition is required for 
$\liminf_{x\to\infty}x^{-1}\log P(X>x)=\sigma_0$. 
In our Theorem $\ref{theo:3}$, the additional condition is the 
analytic property of $\varphi(s)$. 
The example $X^\ast$ in Theorem $\mbox{\rm\ref{theo:4}}$ is in a sense 
pathological, so we can expect that the exponential decay is 
guaranteed by some weak condition. 
\end{remark}

\section{Tauberian Theorems of Laplace Transform}
Theorem \ref{theo:3} deals with an issue of how the analytic 
properties of the Laplace-Stieltjes transform $\varphi(s)$ 
determines the asymptotic behavior of the tail probability 
$P(X>x)$. $\varphi(s)$ converges in the region $\Re s>\sigma_0$ 
and defines an analytic function in this region. 
Thus, according to Widder \cite{wid}, p.40, Theorem 2.2b, 
we see that $P(X>x)=o(e^{\sigma x})$ as $x\to\infty$ for 
$\sigma_0<\sigma<0$. 
The main problem is whether $P(X>x)$ decays as $P(X>x)=O(e^{\sigma_0x})$ 
as $x\to\infty$. 
So, it is appropriate to apply Tauberian theorems of 
Laplace transform to this problem. 

In general, the relation between a function $f$ and its transform 
$Tf$ (such as power series, Laplace transform, etc.) is 
investigated by Abelian theorems or Tauberian theorems. 
In Abelian theorems, the asymptotic behavior of $Tf$ is 
studied from the asymptotic behavior of $f$. Conversely, 
in Tauberian theorems, 
the asymptotic behavior of $f$ is studied from that of 
$Tf$. In Tauberian theorems, generally, some additional 
condition is required for $f$. Such an additional 
condition is called a Tauberian condition. See \cite{kor} 
for the survey of the history and recent developments of 
Tauberian theory.

The following is a well-known Tauberian theorem. 

\medskip

\noindent{\bf Tauberian Theorem} ({\it Widder} \cite{wid}, {\it p.}187, {\it Theorem} 3b) 
{\it For a normalized function $\mu(x)$ of bounded 
variation in $[0,L]$ for every $L>0$, let the integral 
\begin{eqnarray}
\varphi(s)=\int_0^\infty e^{-sx}d\mu(x) \nonumber
\end{eqnarray}
exist for $s>0$ and let $\lim_{s\to0+}\varphi(s)=A$. 
Then 
\begin{eqnarray}
\lim_{x\to\infty}\mu(x)=A \nonumber
\end{eqnarray}
holds if and only if}
\begin{eqnarray}
\ds\int_0^xtd\mu(t)=o(x),\ x\to\infty.\label{eqn:TC}
\end{eqnarray}

\bigskip

In the above theorem, Tauberian condition is (\ref{eqn:TC}). 
So, in this theorem, to study the asymptotic 
property of a function $\mu$, other asymptotic property is 
assumed. I think that this type of Tauberian theorem is 
not easy to apply to some practical problems. 

Ikehara first succeeded \cite{ike} in removing such asymptotic 
conditions from Tauberian theorem, instead, he posed some 
analytic properties of $\varphi(s)$ on the boundary of 
convergence region. 

\subsection{Ikehara's Tauberian Theorem}
The following is Ikehara's Tauberian theorem, in which 
an analytic property of Laplace-Stieltjes transform is 
assumed. The Tauberian condition is the non-decreasing 
property of the function $S(t)$. 

\bigskip

\noindent{\bf Theorem} {\it $($Ikehara $\mbox{\rm\cite{ike}}$, see also 
$\mbox{\rm\cite{kor}})$ Let 
$S(t)$ vanish for $t<0$, be non-decreasing, right continuous, 
and the integral 
\begin{eqnarray}
\varphi(s)=\int_0^\infty e^{-st}dS(t),\ s=\sigma+i\tau
\end{eqnarray}
exist for $\sigma>1$. There exists a constant $A$ such that 
the analytic function 
\begin{eqnarray}
\varphi_0(s)=\varphi(s)-\frac{A}{s-1},\ \Re s>1
\end{eqnarray}
converges as $\sigma\downarrow 1$ to the boundary function 
$\varphi_0(1+i\tau)$ uniformly $($or in $L^1$$)$ for 
$-\lambda<\tau<\lambda$ with any $\lambda>0$. Then we have}
\begin{eqnarray}
\lim_{t\to\infty}e^{-t}S(t)=A. \label{eqn:ikehara}
\end{eqnarray}

\subsection{Finite Form of Ikehara's Theorem by Graham-Vaaler}
In Ikehara's theorem, since $\lambda>0$ is arbitrary, 
$\varphi(s)$ is assumed to be analytic on the whole axis 
of convergence $\Re s=\sigma_0$ except the pole $s=\sigma_0$. 
As mentioned previously, because the axis of convergence 
is not compact, it is difficult to check whether $\varphi(s)$ 
satisfies the theorem assumption or not. The following 
extension by Graham-Vaaler \cite{gra} solves this 
difficulty by relaxing the limit (\ref{eqn:ikehara}) 
of $e^{-t}S(t)$. This theorem is called a finite form 
of Ikehara's theorem because $\lambda$ is restricted 
to some range of values. 

We make preliminary definitions in order to state 
Graham-Vaaler's theorem. 

For $\omega>0$, define a function $E_\omega(t)$ by 
\begin{eqnarray}
E_\omega(t)=\left\{\begin{array}{ll}
e^{-\omega t}, & t\ge 0,\\
0, & t<0.\end{array}\right.\label{eqn:Eomega}
\end{eqnarray}
For $\lambda>0$, a real function $f(x)$ is {\it of type} 
$\lambda$ if $f(x)$ is the restriction to $\mathbb R$ 
of an entire function $f(z)$ of exponential type $\lambda$. 
An entire function $f(z)$ is {\it of exponential type} 
$\lambda$ \cite{rud} if it satisfies 
\begin{eqnarray}
|f(z)|\le C\exp(\lambda|z|),\ z\in\mathbb{C},\ C>0,\ \lambda>0.
\end{eqnarray}

A function $f(x)$ is a {\it majorant} for a function 
$g(x)$ if $f(x)\ge g(x)$ for any $x\in\mathbb R$, and 
$f(x)$ is a {\it minorant} for $g(x)$ if $f(x)\le g(x)$ for 
any $x\in\mathbb R$. 

\bigskip

\noindent{\bf Theorem} {\it$($Graham-Vaaler $\mbox{\rm\cite{gra}}$, 
see also $\mbox{\rm\cite{kor}})$ Let $S(t)$ vanish for $t<0$, be non-decreasing, 
right continuous, and the integral 
\begin{eqnarray}
\varphi(s)=\int_0^\infty e^{-st}dS(t),\ s=\sigma+i\tau
\end{eqnarray}
exist for $\sigma>1$. There exists a constant $A$ 
such that the analytic function 
\begin{eqnarray}
\varphi_0(s)=\varphi(s)-\frac{A}{s-1},\ \Re s>1
\end{eqnarray}
converges as $\sigma\downarrow 1$ to the boundary function 
$\varphi_0(1+i\tau)$ uniformly $($or in $L^1$$)$ for 
$-\lambda<\tau<\lambda$ with some $\lambda>0$. Then, 
for any majorant $M(t)$ for $E_1(t)$ of type $\lambda$ 
and any minorant $m(t)$ for $E_1(t)$ of type $\lambda$, 
we have}
\begin{align}
A\int_{-\infty}^\infty m(t)dt&\le\liminf_{t\to\infty}e^{-t}S(t)\\
&\le\limsup_{t\to\infty}e^{-t}S(t)\le A\int_{-\infty}^\infty M(t)dt.
\end{align}
%

%
%
%
%
%

\section{Main Theorem}
We write below our main theorem again. By this theorem, we 
do not need to check the 
location of singularities of the LS transform $\varphi(s)$ 
and if $\varphi(s)$ is meromorphic then the assumption of 
the theorem is necessarily satisfied.

\medskip

\noindent{\bf Theorem \ref{theo:3}}\ {\it Let $X$ be a non-negative 
random variable, and $F(x)=P(X\le x)$ be the probability 
distribution function of $X$. Let 
\begin{eqnarray}
\varphi(s)=\int_0^\infty e^{-sx}dF(x)
\end{eqnarray}
be the Laplace-Stieltjes transform of $F(x)$ and $\sigma_0$ 
be the abscissa of convergence of $\varphi(s)$. 
We assume $-\infty<\sigma_0<0$. If $s=\sigma_0$ is a pole 
of $\varphi(s)$, then we have}
\begin{eqnarray}
\lim_{x\to\infty}\frac{1}{x}\log P(X>x)=\sigma_0.
\end{eqnarray}

\bigskip

It is possible to give a proof for arbitrary order of poles, but 
the description becomes very complicated, so we will 
prove only for the pole of order 2. A proof for higher order 
poles is easily obtained from the proof for the pole 
of order 2.

\subsection{Preliminary Lemmas for the Proof of Main Theorem}
For $\omega>0$, we define functions
\begin{align}
R(v)&=\frac{v^3}{1-e^{-v}},\ v\in{\mathbb R},\\
\tilde{R}_\omega(v)&=R(v+\omega)-R(\omega)-R'(\omega)v
-\frac{R''(\omega)}{2}v^2,\ v\in{\mathbb R}.\label{eqn:r_tilde}
\end{align}

By calculation, we have 
\begin{align}
R'(v)&=\frac{3v^2}{1-e^{-v}}-\frac{v^3e^{-v}}{(1-e^{-v})^2},
\label{eqn:derivative1}\\
R''(v)&=\frac{6v}{1-e^{-v}}-\frac{6v^2e^{-v}}{(1-e^{-v})^2}
+\frac{2v^3e^{-2v}}{(1-e^{-v})^3}+\frac{v^3e^{-v}}{(1-e^{-v})^2},
\label{eqn:derivative2}\\
R'''(v)&=\frac{6}{1-e^{-v}}-\frac{18ve^{-v}}{(1-e^{-v})^2}
+\frac{18v^2e^{-2v}}{(1-e^{-v})^3}+\frac{9v^2e^{-v}}{(1-e^{-v})^2}
-\frac{6v^3e^{-3v}}{(1-e^{-v})^4}
-\frac{6v^3e^{-2v}}{(1-e^{-v})^3}\nonumber\\
&\ \ -\frac{v^3e^{-v}}{(1-e^{-v})^2},
\label{eqn:derivative3}\\
R^{(4)}(v)&=\frac{-24e^{-v}}{(1-e^{-v})^2}
+\frac{72ve^{-2v}}{(1-e^{-v})^3}
+\frac{36ve^{-v}}{(1-e^{-v})^2}-\frac{72v^2e^{-3v}}{(1-e^{-v})^4}
-\frac{72v^2e^{-2v}}{(1-e^{-v})^3}
-\frac{12v^2e^{-v}}{(1-e^{-v})^2}\nonumber\\
&\ \ +\frac{24v^3e^{-4v}}{(1-e^{-v})^5}
+\frac{36v^3e^{-3v}}{(1-e^{-v})^4}
+\frac{14v^3e^{-2v}}{(1-e^{-v})^3}+\frac{v^3e^{-v}}{(1-e^{-v})^2}.
\label{eqn:derivative4}
\end{align}

We have the following lemmas.

\medskip

\begin{lemma}\label{lem:1}
There exists $\omega_0>0$ such that for any $\omega\ge\omega_0$, 
\begin{eqnarray}\label{eqn:claim:lem1}
\left\{\begin{array}{ll}
\tilde{R}_\omega(v)\ge 0, & v\ge 0,\\[2mm]
\tilde{R}_\omega(v)\le 0, & v<0.
\end{array}\right.
\end{eqnarray}
\end{lemma}

\bigskip

\noindent{\bf Proof:} From (\ref{eqn:derivative3}), 
we have $R'''(v)\simeq 6$ for sufficiently large $v>0$, 
or more precisely, for any $\epsilon>0$ there exists 
$v_0>0$ such that $|R'''(v)-6|<\epsilon$ holds for 
$v\geq v_0$. Hence, by taking $\epsilon$ sufficiently 
small, we have $R'''(v)>5>0$ for 
$v\geq v_0$, especially, $R''(v)$ is monotonic increasing 
for $v\geq v_0$. Further, from (\ref{eqn:derivative2}) 
we have 
\begin{align}
\lim_{v\to\infty}R''(v)&=\infty,\label{eqn:lem1:1}\\
\lim_{v\to-\infty}R''(v)&=0.\label{eqn:lem1:2}
\end{align}
Write $C\equiv\max_{v\leq v_0}R''(v)$, which is finite 
by (\ref{eqn:lem1:2}). From (\ref{eqn:lem1:1}), there 
exists $\omega_0$ with $R''(\omega_0)=C+1$ and 
$\omega_0\geq v_0$. Then, we have 
\begin{align}
\max_{v\leq\omega_0}R''(v)&=\max\left(\max_{v\leq v_0}R''(v), 
\max_{v_0\leq v\leq\omega_0}R''(v)\right)\\
&=\max(C, C+1)\\
&=C+1\\
&=R''(\omega_0).
\end{align}
Since $R''(v)$ is monotonic increasing for 
$v\geq\omega_0(\geq v_0)$, we have
\begin{align}
\max_{v\leq\omega}R''(v)=R''(\omega),\ 
\forall\omega\geq\omega_0.\label{eqn:lem1:3}
\end{align}
Next, for $\omega\geq\omega_0$, we have
\begin{align}
\max_{v\geq\omega}R''(v)=R''(\omega),\ 
\forall\omega\geq\omega_0\label{eqn:lem1:4}
\end{align}
because $R''(v)$ is monotonic increasing for $v\geq\omega
(\geq\omega_0\geq v_0)$. From (\ref{eqn:r_tilde}), we have 
$\tilde{R}''_\omega(v)=R''(v+\omega)-R''(\omega)$, 
thus for $\omega\geq\omega_0$, we have from (\ref{eqn:lem1:3}), 
(\ref{eqn:lem1:4}),
\begin{align}
\tilde{R}''_\omega(v)
\left\{\begin{array}{ll}
\ge 0, & v\ge 0,\\
\le 0, & v<0.
\end{array}\right.
\end{align}
From $\tilde{R}'_\omega(0)=0$, we have 
$\tilde{R}'_\omega(v)\ge0,\ v\in{\mathbb R}$. 
So, $\tilde{R}_\omega(v)$ is monotonic increasing in $v$. 
Then, $\tilde{R}_\omega(0)=0$ implies (\ref{eqn:claim:lem1}). 

\begin{lemma}\label{lem:2}
For sufficiently large $v>0$,
\begin{align}
R'(v)&<\frac{3v^2}{1-e^{-v}},\label{eqn:lem2:1}\\
R'''(v)&<\frac{6}{1-e^{-v}}.\label{eqn:lem2:2}
\end{align}
\end{lemma}

\bigskip

\noindent{\bf Proof:} (\ref{eqn:lem2:1}) is easily obtained 
from (\ref{eqn:derivative1}). We have from (\ref{eqn:derivative3}) 
\begin{align}
R'''(v)&=\ds\frac{6}{1-e^{-v}}-v^3e^{-v}\left(1+O(v^{-1})\right),\ v\to\infty,
\end{align}
thus, (\ref{eqn:lem2:2}) holds.

\bigskip

\begin{lemma}\label{lem:3}
There exists $\omega_0>0$ such that for any $\omega\ge\omega_0$, 
\begin{eqnarray}
\tilde{R}_\omega(v)\le\frac{1}{1-e^{-\omega}}v^3,\ v\ge 0.
\end{eqnarray}
\end{lemma}

\medskip

\noindent{\bf Proof:} By the mean value theorem and Lemma \ref{lem:2}, 
there exists $v_0$ with $\omega\le v_0\le v+\omega$ such that 
\begin{align}
\frac{\tilde{R}_\omega(v)}{v^3}&=\frac{R'''(v_0)}{3!}\nonumber\\
&<\frac{1}{3!}\cdot\frac{6}{1-e^{-v_0}}\\
&\le\ds\frac{1}{1-e^{-\omega}}.
\end{align}

\bigskip

\begin{lemma}\label{lem:4}
There exists $\omega_0>0$ such that for any $\omega\ge\omega_0$, 
\begin{eqnarray}
-\int_{-\infty}^0\tilde{R}_\omega(v)e^{tv}dv
<\frac{R'''(\omega)}{t^4},\ t>0.
\label{eqn:claim:lem4}
\end{eqnarray}
\end{lemma}

\bigskip

\noindent{\bf Proof:} By calculation, for $\omega>0$ 
\begin{align}
\int_{-\infty}^0\tilde{R}_\omega(v)e^{tv}dv
=\int_{-\infty}^0R(v+\omega)e^{tv}dv-\frac{R(\omega)}{t}
+\frac{R'(\omega)}{t^2}-\frac{R''(\omega)}{t^3},\label{eqn:lem4:1}
\end{align}
and from integration by parts 
\begin{align}
\int_{-\infty}^0R(v+\omega)e^{tv}dv
=\frac{R(\omega)}{t}-\frac{R'(\omega)}{t^2}+\frac{R''(\omega)}{t^3}
-\frac{R'''(\omega)}{t^4}
+\frac{1}{t^4}\int_{-\infty}^0
R^{(4)}(v+\omega)e^{tv}dv,\label{eqn:lem4:2}
\end{align}
where $R^{(4)}(v)$ denotes the fourth derivative of $R(v)$. 
Hence from (\ref{eqn:lem4:1}), (\ref{eqn:lem4:2}), we have 
\begin{align}
\frac{R'''(\omega)}{t^4}+\int_{-\infty}^0\tilde{R}_\omega(v)e^{tv}dv
&=\frac{1}{t^4}\int_{-\infty}^0R^{(4)}(v+\omega)e^{tv}dv\nonumber\\
&=\frac{e^{-\omega t}}{t^4}\int_{-\infty}^\omega R^{(4)}(v)e^{tv}dv.\nonumber
\end{align}
Then, we will prove that there exists $\omega_0>0$ such that 
for any $\omega\ge\omega_0$, 
\begin{eqnarray}
\int_{-\infty}^\omega R^{(4)}(v)e^{tv}dv>0,\ t\ge0.\label{eqn:lem4:3}
\end{eqnarray}
We divide into two cases $t>1$ and $0\le t\le1$. 

First, let us consider the case $t>1$. From (\ref{eqn:derivative4}), 
we have 
\begin{align}
R^{(4)}(v)=\left\{\begin{array}{ll}
v^3e^{-v}\left(1+O(v^{-1})\right), & v\to\infty,\\[1mm]
O(v^3e^v), & v\to -\infty, \label{eqn:lem4:4}
\end{array}\right.
\end{align}
especially, $R^{(4)}(v)$ is bounded in ${\mathbb R}$. 
From (\ref{eqn:lem4:4}) we see that there exist constants 
$C>0$, $v_0>0$ with
\begin{align}
R^{(4)}(v)\ge\left\{\begin{array}{ll}
\ds\frac{1}{2}v^3e^{-v}, & v\ge v_0,\\[2mm]
-C, & v<v_0.
\end{array}\right.
\end{align}
Thus, for $t>1$, 
\begin{align}
\int_{-\infty}^\omega R^{(4)}(v)e^{tv}dv
&\ge -C\int_{-\infty}^{v_0}e^{tv}dv
+\frac{1}{2}\int_{v_0}^\omega v^3e^{(t-1)v}dv\\[2mm]
&\ge -Ce^{tv_0}
+\frac{v_0^3}{2}\frac{e^{\omega(t-1)}-e^{v_0(t-1)}}{t-1}\\[2mm]
&=\frac{e^{\omega(t-1)}-e^{v_0(t-1)}}{t-1}\left\{\frac{v_0^3}{2}
-\frac{Ce^{tv_0}}{e^{\omega(t-1)}-e^{v_0(t-1)}}\right\}.\label{eqn:lem4:4-2}
\end{align}
For any $\epsilon>0$, putting $\omega_0=v_0+e^{v_0}/\epsilon$, 
we have for any $\omega\ge\omega_0$, 
\begin{align}
e^{\omega(t-1)} &\ge e^{(v_0+e^{v_0}/\epsilon)(t-1)} \\
&\ge e^{v_0(t-1)}\left\{1+\frac{e^{v_0}}{\epsilon}(t-1)\right\},\ t>1,
\label{eqn:lem4:5}
\end{align}
or 
\begin{align}
\frac{(t-1)e^{v_0t}}{e^{\omega(t-1)}-e^{v_0(t-1)}}\le\epsilon,\ t>1.
\label{eqn:lem4:4-3}
\end{align}
Therefore, we have from (\ref{eqn:lem4:4-2}),(\ref{eqn:lem4:4-3}), 
by 
taking $\epsilon$ sufficiently small, for $\omega\geq\omega_0$, 
\begin{align}
\int_{-\infty}^\omega R^{(4)}(v)e^{tv}dv>0,\ t>1. \label{eqn:lem4:6}
\end{align}

Next, we consider $0\le t\le 1$. We first restrict $0<t<1$. 
The cases $t=0$ and $t=1$ will be considered later. 
We calculate
\begin{eqnarray}
\int_{-\infty}^\infty R^{(4)}(v)e^{tv}dv,\ 0<t<1.
\end{eqnarray}
Here, we put $t=-s$ and 
\begin{eqnarray}
\phi(s)=\int_{-\infty}^\infty R^{(4)}(v)e^{-sv}dv,\ -1<s<0.
\end{eqnarray}
$\phi(s)$ is the bilateral Laplace transform \cite{wid} of $R^{(4)}(v)$. 
The region of convergence is $-1<\Re s<1$. 

We apply the following theorem. 

\medskip

\noindent{\bf Theorem A} {\it$($Widder $\mbox{\rm\cite{wid}}$, p.$239$, 
Theorem $3$d$)$ Let $\alpha(v)$ be of bounded variation 
on any finite interval. If the integral
\begin{eqnarray}
f(s)=\int_{-\infty}^\infty e^{-sv}d\alpha(v)
\end{eqnarray}
exists for $s=s_0,\ \Re s_0<0$, and $\alpha(\infty)=0$, 
then}
\begin{eqnarray}
f(s_0)=s_0\ds\int_{-\infty}^\infty e^{-s_0v}\alpha(v)dv.
\end{eqnarray}
\noindent(end of theorem citation)

\bigskip

Let us define a step function $\Delta(v)$ as
\begin{eqnarray}
\Delta(v)=\left\{\begin{array}{ll}
1, & v\ge 0,\\
0, & v<0,
\end{array}\right.
\end{eqnarray}
and define $\alpha_1(v)=R'''(v)-6\Delta(v)$. Then, 
\begin{eqnarray}
\int_{-\infty}^\infty e^{-sv}d\alpha_1(v)=\int_{-\infty}^\infty 
e^{-sv}R^{(4)}(v)dv-6<\infty, \label{eqn:lem4:7}
\end{eqnarray}
and $\alpha_1(\infty)=0$, so $\alpha_1(v)$ satisfies the 
assumptions of the above Theorem A. Hence, we have 
from Theorem A and (\ref{eqn:lem4:7}) 
\begin{eqnarray}
\int_{-\infty}^\infty e^{-sv}R^{(4)}(v)dv-6=s\int_{-\infty}^\infty e^{-sv}
\alpha_1(v)dv. \label{eqn:lem4:8}
\end{eqnarray}

Next, define $\alpha_2(v)=R''(v)-6v\Delta(v)$, then 
from (\ref{eqn:lem4:7}), (\ref{eqn:lem4:8})
\begin{eqnarray}
\int_{-\infty}^\infty e^{-sv}d\alpha_2(v)=\int_{-\infty}^\infty e^{-sv}
\alpha_1(v)dv<\infty, \label{eqn:lem4:9}
\end{eqnarray}
and $\alpha_2(\infty)=0$, so $\alpha_2(v)$ satisfies the 
assumptions of Theorem A. Hence, we have from Theorem A and 
(\ref{eqn:lem4:9})
\begin{eqnarray}
\int_{-\infty}^\infty e^{-sv}\alpha_1(v)dv=s\int_{-\infty}^\infty e^{-sv}
\alpha_2(v)dv. \label{eqn:lem4:10}
\end{eqnarray}

In a similar way, by defining $\alpha_3(v)=R'(v)-3v^2\Delta(v)$, 
$\alpha_4(v)=R(v)-v^3\Delta(v)$, we have from Theorem A
\begin{align}
\int_{-\infty}^\infty e^{-sv}\alpha_2(v)dv&=s\int_{-\infty}^\infty e^{-sv}
\alpha_3(v)dv, \label{eqn:lem4:11}\\
\int_{-\infty}^\infty e^{-sv}\alpha_3(v)dv&=s\int_{-\infty}^\infty e^{-sv}
\alpha_4(v)dv. \label{eqn:lem4:12}
\end{align}
Therefore, from (\ref{eqn:lem4:8}), (\ref{eqn:lem4:10}), 
(\ref{eqn:lem4:11}), (\ref{eqn:lem4:12}) we have
\begin{align}
\int_{-\infty}^\infty e^{-sv}R^{(4)}(v)dv-6
&=s^4\int_{-\infty}^\infty e^{-sv}\alpha_4(v)dv\nonumber\\
&=s^4\left\{\int_{-\infty}^0\frac{v^3e^{-sv}}{1-e^{-v}}dv
+\int_0^\infty\frac{v^3e^{-sv}}{e^v-1}dv\right\},\ -1<s<0.\nonumber
\end{align}
By the change of variables $t=-s$, 
\begin{align}
\int_{-\infty}^\infty R^{(4)}(v)e^{tv}dv
&=6+t^4\left\{\int_{-\infty}^0\frac{v^3e^{tv}}{1-e^{-v}}dv
+\int_0^\infty\frac{v^3e^{tv}}{e^v-1}dv\right\}\\
&=6+t^4\left\{\int_0^\infty\frac{v^3e^{-tv}}{e^v-1}dv
+\int_0^\infty\frac{v^3e^{tv}}{e^v-1}dv\right\}\\
&>6,\ 0<t<1.\label{eqn:lem4:13}
\end{align}

Next consider $t=0$. Define a function $g(v)$ as 
\begin{eqnarray}
g(v)=\left\{\begin{array}{ll}
\left|R^{(4)}(v)\right|e^{v/2}, & v\ge0,\\[3mm]
\left|R^{(4)}(v)\right|, & v<0.
\end{array}\right.
\end{eqnarray}
Then from (\ref{eqn:lem4:4}), $g(v)$ is integrable and 
$|R^{(4)}(v)|e^{tv}\le g(v),\ v\in{\mathbb R}$, for $0<t<1/2$.
By the dominating convergence theorem 
\begin{align}
\int_{-\infty}^\infty R^{(4)}(v)dv
&=\lim_{t\to0+}\int_{-\infty}^\infty R^{(4)}(v)e^{tv}dv\nonumber\\
&=6+\lim_{t\to0}t^4\cdot\lim_{t\to0+}\left(\int_{-\infty}^0
\ds\frac{v^3e^{tv}}{1-e^{-v}}dv
+\int_0^\infty \ds\frac{v^3e^{tv}}{e^v-1}dv\right)\nonumber\\
&=6+\lim_{t\to0}t^4\cdot\left(\int_{-\infty}^0\ds\frac{v^3}{1-e^{-v}}dv
+\int_0^\infty \ds\frac{v^3}{e^v-1}dv\right)\nonumber\\
&=6. \label{eqn:lem4:14}
\end{align}

Last, for $t=1$, we have from (\ref{eqn:lem4:4})
\begin{eqnarray}
\int_{-\infty}^\infty R^{(4)}(v)e^vdv=\infty.\label{eqn:lem4:15}
\end{eqnarray}
Summarizing (\ref{eqn:lem4:13}), (\ref{eqn:lem4:14}), 
(\ref{eqn:lem4:15}), we obtain 
\begin{eqnarray}
\int_{-\infty}^\infty R^{(4)}(v)e^{tv}dv\ge 6,\ 0\le t\le1.
\label{eqn:lem4:16}
\end{eqnarray}
Now we show that there exists $\omega_0>0$ such that 
for any $\omega\ge\omega_0$, 
\begin{eqnarray}
\int_{-\infty}^\omega R^{(4)}(v)e^{tv}dv>5,\ 0\le t\le1.
\label{eqn:lem4:b}
\end{eqnarray}
Define 
\begin{eqnarray}
\phi_\omega(t)=\int_{-\infty}^\omega R^{(4)}(v)e^{tv}dv,
\ 0\le t\le 1,\ \omega>0.\label{eqn:phi_omega}
\end{eqnarray}

For a fixed $\omega>0$, the integral in (\ref{eqn:phi_omega}) 
converges for $t>-1$. In fact, from (\ref{eqn:derivative4}), 
for any $\delta>0$, 
\begin{align}
\left|R^{(4)}(v)\right|\leq e^{(1-\delta)v}\label{eqn:lem4:a}
\end{align}
holds for any $v<0$ with sufficiently large absolute value. 
By the change of variables $v=\omega-u$ in the integral 
in (\ref{eqn:phi_omega}), we have 
\begin{align}
\phi_\omega(t)=e^{t\omega}\int_0^\infty R^{(4)}(\omega-u)
e^{-tu}du.\label{eqn:lem4:c}
\end{align}
Hence from (\ref{eqn:lem4:a}), the integral in (\ref{eqn:lem4:c}) 
converges for $t>-1$. Then from Widder \cite{wid}, p.57, Theorem 5a, 
we see that $\phi_\omega(t)$ is analytic for $\Re t>-1$, 
especially $\phi_\omega(t)$ is continuous for $0\leq t\leq 1$.

For a fixed $t$ with $0\le t\le1$, from (\ref{eqn:lem4:4}) 
there exists $\omega^\ast>0$ (which does not depend on $t$) 
such that for $\omega>\omega^\ast$, $\phi_\omega(t)$ is 
monotonic increasing in $\omega$. 

From (\ref{eqn:lem4:16}) for any $t$ with $0\le t\le 1$, 
there exists $\omega(t)\ge\omega^\ast$ such that 
\begin{eqnarray}
\phi_{\omega(t)}(t)=\int_{-\infty}^{\omega(t)}R^{(4)}(v)e^{tv}dv>5.
\end{eqnarray}
By the continuity property of $\phi_{\omega(t)}(u)$ in $u$, 
there is $\epsilon=\epsilon(t, \omega(t))>0$ such that 
\begin{eqnarray}
\phi_{\omega(t)}(u)>5,\ u\in(t-\epsilon,\ t+\epsilon)\cap[0,1].
\end{eqnarray}
Define $U_t=(t-\epsilon,\ t+\epsilon)\cap[0,1]$, then 
$U_t$ is an open set in $[0,1]$ and $t\in U_t$, so 
$\bigcup_{\,0\le t\le1}U_t=[0,1]$. Since $[0,1]$ is 
compact, there exist $t_1,\cdots,t_n\in[0,1]$ with 
$U_{t_1}\cup\cdots\cup U_{t_n}=[0,1]$. 
Let $\omega_0=\max_{1\le k\le n}\omega(t_k)$, then for any 
$\omega\ge\omega_0$, 
\begin{align}
\phi_\omega(t) & \ge \phi_{\omega(t_k)}(t),\ t\in U_{t_k}\\
&>5 \label{eqn:lem4:17}
\end{align}
holds for $t\in[0,1]$. From (\ref{eqn:lem4:6}), (\ref{eqn:lem4:17}) 
we have completed the proof of Lemma \ref{lem:4}.

\bigskip

Now, for $\omega>0$, define two functions $M_\omega(t)$, 
$m_\omega(t)$ as follows.
\begin{align}
M_\omega(t)&=\frac{1}{6}\left(\frac{\sin\pi t}{\pi}\right)^4Q_\omega(t),
\ -\infty<t<\infty, \label{eqn:lem5pre1}\\
Q_\omega(t)&=\sum_{n=0}^\infty e^{-n\omega}\left\{\frac{6}{(t-n)^4}
-\frac{6\omega}{(t-n)^3}+\frac{3\omega^2}{(t-n)^2}
-\frac{\omega^3}{t-n}\right\}
+\frac{R(\omega)}{t}-\frac{R'(\omega)}{t^2}+\frac{R''(\omega)}{t^3},\\
m_\omega(t)&=M_\omega(t)-\ds\frac{1}{1-e^{-\omega}}
\left(\frac{\sin\pi t}{\pi t}\right)^4,\ -\infty<t<\infty.\label{eqn:lem5pre2}
\end{align}

\begin{lemma}\label{lem:5}
There exists $\omega_0>0$ such that for any $\omega\ge\omega_0$, 
$M_\omega(t)$ is a majorant for $E_\omega(t)$ of 
type $4\pi$ and $\int_{-\infty}^\infty M_\omega(t)dt<\infty$, 
moreover, $m_\omega(t)$ is a minorant for $E_\omega(t)$ 
of type $4\pi$ and $\int_{-\infty}^\infty m_\omega(t)dt>0$.
\end{lemma}

\bigskip

\noindent{\bf Proof:} First, consider $t<0$. We have 
\begin{align}
\int_0^\infty R(v+\omega)e^{tv}dv&=\int_\omega^\infty 
R(u)e^{t(u-\omega)}du\\
&=e^{-\omega t}\int_\omega^\infty u^3e^{tu}\sum_{n=0}^\infty 
e^{-nu}du.
\end{align}
From the monotone convergence theorem and integration by 
parts 
\begin{align}
\int_0^\infty R(v+\omega)e^{tv}dv
=\sum_{n=0}^\infty e^{-n\omega}
\left\{\frac{6}{(t-n)^4}-\frac{6\omega}{(t-n)^3}
+\frac{3\omega^2}{(t-n)^2}-\frac{\omega^3}{t-n}\right\}.\nonumber
\end{align}
Then, we have from (\ref{eqn:r_tilde})
\begin{align}
\int_0^\infty \tilde{R}_\omega(v)e^{tv}dv=Q_\omega(t),\ t<0.
\label{eqn:lem5:1}
\end{align}
From Lemmas \ref{lem:1}, \ref{lem:3}, for any $\omega\ge\omega_0$ 
\begin{align}
0\le\tilde{R}_\omega(v)\le\frac{1}{1-e^{-\omega}}v^3,\ v\ge0,
\end{align}
hence from (\ref{eqn:lem5:1}) 
\begin{align}
0\le Q_\omega(t)\le\frac{6}{1-e^{-\omega}}\cdot\frac{1}{t^4},\ t<0,
\end{align}
and by (\ref{eqn:lem5pre1}) we have 
\begin{align}
0\le M_\omega(t)\le\frac{1}{1-e^{-\omega}}
\left(\frac{\sin\pi t}{\pi t}\right)^4,\ t<0. \label{eqn:lem5:2}
\end{align}

Next, consider $t>0,\ t\notin\mathbb{Z}$. Let us calculate 
$-\int_{-\infty}^0\tilde{R}_\omega(v)e^{tv}dv$. We have 
\begin{eqnarray}
\int_{-\infty}^0 R(v+\omega)e^{tv}dv=
e^{-\omega t}\left\{\int_{-\infty}^0+\int_0^\omega\right\}
R(u)e^{tu}du.
\end{eqnarray}
From the monotone convergence theorem and integration by 
parts 
\begin{eqnarray}
\int_{-\infty}^0 R(u)e^{tu}du=\sum_{n=1}^\infty \frac{6}{(t+n)^4},
\end{eqnarray}
and 
\begin{align}
\int_0^\omega R(u)e^{tu}du=\sum_{n=0}^\infty\frac{6}{(t-n)^4}
+e^{\omega t}\sum_{n=0}^\infty e^{-n\omega}\left\{\frac{\omega^3}{t-n}
-\frac{3\omega^2}{(t-n)^2}+\frac{6\omega}{(t-n)^3}-\frac{6}{(t-n)^4}
\right\}.\nonumber
\end{align}

We obtain a formula
\begin{align}
\sum_{n=-\infty}^\infty\frac{1}{(t-n)^4}
=\left(\frac{\pi}{\sin\pi t}\right)^4\label{eqn:ahlfors_formula}
\end{align}
by the same argument as in Ahlfors \cite{ahl}, p.188. From 
(\ref{eqn:ahlfors_formula}),
\begin{align}
-\int_{-\infty}^0\tilde{R}_\omega(v)e^{tv}dv
=Q_\omega(t)-6e^{-\omega t}\left(\frac{\pi}{\sin\pi t}\right)^4.
\end{align}
Hence, from Lemmas \ref{lem:1}, \ref{lem:4}, for $\omega\ge\omega_0$ 
\begin{align}
0\le Q_\omega(t)-6e^{-\omega t}\left(\frac{\pi}{\sin\pi t}\right)^4
\le\frac{R'''(\omega)}{t^4},
\end{align}
and from Lemma \ref{lem:2}
\begin{align}
e^{-\omega t}\le M_\omega(t)\le e^{-\omega t}
+\frac{1}{1-e^{-\omega}}\left(\frac{\sin\pi t}{\pi t}\right)^4,
\ t>0,\ t\notin{\mathbb Z}.\label{eqn:lem5:3}
\end{align}
By continuity, (\ref{eqn:lem5:3}) holds for all $t>0$.

For $t=0$, $E_\omega(0)=M_\omega(0)=1$. Therefore from 
(\ref{eqn:lem5:2}), (\ref{eqn:lem5:3}) we see that $M_\omega(t)$ 
is a majorant for $E_\omega(t)$ for $\omega\ge\omega_0$.

$M_\omega(t)$ is of type $4\pi$ by the following reason. Consider 
a rectangle $\Gamma_n$ in the complex plane with vertices $(n+1/2)+i$, 
$(-n-1/2)+i$, $(-n-1/2)-i$, $(n+1/2)-i,\ n=1,2,\cdots$. From 
$|\sin\pi z|^2=\sin^2\pi t+\sinh^2\pi y,\ z=t+iy$, we see that $|\sin\pi z|^2$ is 
bounded by a constant on $\Gamma_n$, and $|Q_\omega(z)|$ is 
also bounded by a constant on $\Gamma_n$. Both constants 
are independent on $n$. Thus, by letting $n\to\infty$, 
from the maximum principle, we have for some constant $C$ 
\begin{align}
|M_\omega(z)|\leq C,\ \mbox{\rm for}\ |y|\leq 1.\label{eqn:type1}
\end{align}
Next, since $|Q_\omega(z)|$ is bounded in $|y|\geq 1$, we have 
\begin{align}
|M_\omega(z)|\leq C|\sin\pi z|^4,\ \mbox{\rm for}\ |y|\geq 1.\label{eqn:type2}
\end{align}
By the formula $\sin\pi z=(e^{i\pi z}-e^{-i\pi z})/2i$, we 
have 
\begin{align}
|\sin\pi z|^4\leq e^{4\pi|z|},\ z\in{\mathbb C}.\label{eqn:type3}
\end{align}
From (\ref{eqn:type1}),(\ref{eqn:type2}),(\ref{eqn:type3}), we 
have 
\begin{align}
|M_\omega(z)|\leq Ce^{4\pi|z|},\ z\in{\mathbb C},
\end{align}
which implies $M_\omega(t)$ is of type $4\pi$.

$\int_{-\infty}^\infty M_\omega(t)dt<\infty$ is trivial 
from (\ref{eqn:lem5:2}), (\ref{eqn:lem5:3}), but we 
calculate the integral. Note that 
\begin{align}
\int_{-\infty}^\infty\left\{\frac{R(\omega)}{t}-\sum_{n=0}^\infty 
e^{-n\omega}\frac{\omega^3}{t-n}\right\}
\left(\frac{\sin\pi t}{\pi}\right)^4dt
=-\omega^3\int_{-\infty}^\infty \sum_{n=0}^\infty ne^{-n\omega}
\frac{1}{t(t-n)}\left(\frac{\sin\pi t}{\pi}\right)^4dt.\label{eqn:lem5:4}
\end{align}

Is is easy to check the following inequalities;
\begin{eqnarray*}
\left|t(t-n)\right|\geq
\left\{\begin{array}{ll}
(t-n)^2, & t\geq n,\\
t^2, & t<0,
\end{array}
\right.
\end{eqnarray*}
and 
\begin{align*}
\frac{1}{\left|t(t-n)\right|}\left(\frac{\sin\pi t}{\pi}\right)^4
&=\frac{\sin^2\pi t}{\pi^2}\left|\frac{\sin\pi t}{\pi t}\right|
\left|\frac{\sin\pi(t-n)}{\pi(t-n)}\right|\\
&\leq\frac{1}{\pi^2},\ t\in{\mathbb R}.
\end{align*}
Hence, we have
\begin{align}
\int_{-\infty}^\infty\left|\frac{1}{t(t-n)}\right|
\left(\frac{\sin\pi t}{\pi}\right)^4dt
&\le \int_{-\infty}^0\frac{1}{t^2}\left(\frac{\sin\pi t}{\pi}\right)^4dt
+\int_0^n \frac{1}{\pi^2}dt+\int_n^\infty\frac{1}{(t-n)^2}
\left(\frac{\sin\pi t}{\pi}\right)^4dt\nonumber\\
&= \int_{-\infty}^\infty\frac{1}{t^2}\left(\frac{\sin\pi t}{\pi}\right)^4dt
+\frac{n}{\pi^2},\nonumber
\end{align}
From the dominating convergence theorem, 
we can calculate (\ref{eqn:lem5:4}) as follows 
(by considering Cauchy's principal value)
\begin{align}
&\int_{-\infty}^\infty\left\{\frac{R(\omega)}{t}-\sum_{n=0}^\infty 
e^{-n\omega}\frac{\omega^3}{t-n}\right\}
\left(\frac{\sin\pi t}{\pi}\right)^4dt\nonumber\\
&=\omega^3\sum_{n=0}^\infty e^{-n\omega}\left\{\int_{-\infty}^\infty
\frac{1}{t}\left(\frac{\sin\pi t}{\pi}\right)^4dt
-\int_{-\infty}^\infty\frac{1}{t-n}\left(\frac{\sin\pi t}{\pi}\right)^4dt
\right\}\\
&=0.
\end{align}
Therefore, we have
\begin{align}
\int_{-\infty}^\infty M_\omega(t)dt
&=\frac{1}{6}\sum_{n=0}^\infty e^{-n\omega}\int_{-\infty}^\infty 
\left(\frac{\sin\pi t}{\pi}\right)^4\frac{6}{(t-n)^4}dt
+\frac{1}{6}\sum_{n=0}^\infty e^{-n\omega}\int_{-\infty}^\infty 
\left(\frac{\sin\pi t}{\pi}\right)^4\frac{3\omega^2}{(t-n)^2}dt\nonumber\\
&\hspace{10mm}-\frac{1}{6}\int_{-\infty}^\infty
\left(\frac{\sin\pi t}{\pi}\right)^4\frac{R'(\omega)}{t^2}dt\\
&=\frac{1}{1-e^{-\omega}}\int_{-\infty}^\infty
\left(\frac{\sin\pi t}{\pi t}\right)^4dt
+\frac{1}{6\pi^2}\left\{\frac{3\omega^2}{1-e^{-\omega}}
-R'(\omega)\right\}\int_{-\infty}^\infty\frac{\sin^4\pi t}{(\pi t)^2}dt\\
&<\infty.
\end{align}

Next, we prove for $m_\omega(t)$. From (\ref{eqn:lem5pre2}), 
(\ref{eqn:lem5:3}), for $\omega\ge\omega_0$ 
\begin{align}
e^{-\omega t}-\frac{1}{1-e^{-\omega}}\left(\frac{\sin\pi t}{\pi t}\right)^4
\le m_\omega(t)\le e^{-\omega t},\ t\ge0
\end{align}
and from (\ref{eqn:lem5pre2}), (\ref{eqn:lem5:2}),
\begin{align}
-\frac{1}{1-e^{-\omega}}\left(\frac{\sin\pi t}{\pi t}\right)^4
\le m_\omega(t)\le 0,\ t<0
\end{align}
so, $m_\omega(t)$ is a minorant for $E_\omega(t)$ and is 
of type $4\pi$. Since 
\begin{align}
\int_{-\infty}^\infty m_\omega(t)dt
&=\int_{-\infty}^\infty M_\omega(t)dt-\frac{1}{1-e^{-\omega}}
\int_{-\infty}^\infty\left(\frac{\sin\pi t}{\pi t}\right)^4dt\nonumber\\
&=\frac{1}{6\pi^2}\left\{\frac{3\omega^2}{1-e^{-\omega}}-R'(\omega)\right\}
\int_{-\infty}^\infty\frac{\sin^4\pi t}{(\pi t)^2}dt,\nonumber
\end{align}
from Lemma \ref{lem:2} for $\omega\ge\omega_0$
\begin{eqnarray}
\int_{-\infty}^\infty m_\omega(t)dt>0. 
\end{eqnarray}

\subsection{Extension of Graham-Vaaler's Theorem}
Substituting $t\to-\ds\frac{\sigma_0}{\omega}t$ into 
$E_\omega(t)$ defined by (\ref{eqn:Eomega}), we have
\begin{eqnarray}
E_\omega(-\frac{\sigma_0}{\omega}t)=E_{\sigma_0}(t).
\end{eqnarray}
Define $\lambda= 2\pi/\omega$ and 
\begin{align}
M_{\lambda, \sigma_0}(t)&= M_\omega(-\frac{\sigma_0}{\omega}t),\ 
t\in{\mathbb R},\ \omega>0, \\
m_{\lambda, \sigma_0}(t)&= m_\omega(-\frac{\sigma_0}{\omega}t),\ 
t\in{\mathbb R},\ \omega>0,
\end{align}
then both $M_{\lambda, \sigma_0}(t)$ and $m_{\lambda, \sigma_0}(t)$ 
are of type $-2\sigma_0\lambda$. From Lemma \ref{lem:5}, 
for sufficiently small $\lambda>0$, $M_{\lambda, \sigma_0}(t)$ 
and $m_{\lambda, \sigma_0}(t)$ are majorant and minorant 
for $E_{\sigma_0}(t)$, respectively. 

Based on Lemma \ref{lem:5}, we can calculate
\begin{align}
\int_{-\infty}^\infty M_{\lambda, \sigma_0}(t)dt
&=-\frac{\omega}{\pi\sigma_0(1-e^{-\omega})}
\int_{-\infty}^\infty\left(\frac{\sin t}{t}\right)^4dt
-\frac{\omega}{6\sigma_0\pi^3}\left\{\frac{3\omega^2}{1-e^{-\omega}}
-R'(\omega)\right\}\int_{-\infty}^\infty\frac{\sin^4t}{t^2}dt\\
&<\infty,\\[3mm]
\int_{-\infty}^\infty m_{\lambda, \sigma_0}(t)dt
&=-\frac{\omega}{6\sigma_0\pi^3}\left\{\frac{3\omega^2}{1-e^{-\omega}}
-R'(\omega)\right\}\int_{-\infty}^\infty\frac{\sin^4t}{t^2}dt\\
&>0.
\end{align}

We have the following theorem. The proof is an extension 
of the proof of Graham-Vaaler's Tauberian theorem to the 
case that the pole is of order $2$.

\medskip

\begin{theorem}\label{theo:6} Let $X$ be a non-negative 
random variable and $\varphi(s)$ be the Laplace-Stieltjes 
transform of the probability distribution function 
$F(x)$ of $X$. The abscissa of convergence of $\varphi(s)$ 
is denoted by $\sigma_0$ and $-\infty<\sigma_0<0$ is 
assumed. Let $s=\sigma_0$ be a pole of $\varphi(s)$ of 
order $2$ and $\varphi(s)$ be analytic on the interval 
$\{s=\sigma_0+i\tau\,|\,2\sigma_0\lambda<\tau<-2\sigma_0\lambda\}$ 
for some $\lambda>0$ except $s=\sigma_0$. Write $A_2$ as 
the coefficient of $(s-\sigma_0)^{-2}$ in the Laurent 
expansion of $\varphi(s)$ at $s=\sigma_0$. 
Then, we have
\begin{align}
A_2\int_{-\infty}^\infty m_{\lambda, \sigma_0}(t)dt
&\le\liminf_{x\to\infty}x^{-1}e^{-\sigma_0x}P(X>x)\nonumber\\
&\le\limsup_{x\to\infty}x^{-1}e^{-\sigma_0x}P(X>x)\nonumber\\
&\le A_2\int_{-\infty}^\infty M_{\lambda, \sigma_0}(t)dt.\nonumber
\end{align}
\end{theorem}

\medskip

\noindent{\bf Proof:} For the evaluation of 
$P(X>x)=\int_x^\infty dF(t)$, we first evaluate the 
following integral. For $\sigma>\sigma_0$, we have
\begin{align}
e^{-\sigma_0 x}\int_x^\infty e^{-(\sigma-\sigma_0)t}dF(t)
&=\int_x^\infty e^{\sigma_0(t-x)}e^{-\sigma t}dF(t)\nonumber\\
&=\int_0^\infty E_{\sigma_0}(t-x)e^{-\sigma t}dF(t)\label{eqn:second}\\
&\le \int_0^\infty M_{\lambda,\sigma_0}(t-x)e^{-\sigma t}dF(t),\nonumber
\end{align}
if $F(t)$ is continuous at $x$. If $x$ is a point of 
discontinuity of $F(t)$, then the integral in (\ref{eqn:second}) 
does not exist. Since $F(t)$ is right continuous, the 
set of points of discontinuity is at most countable. 
We see that the integral 
$\int_0^\infty M_{\lambda,\sigma_0}(t-x)e^{-\sigma t}dF(t)$ 
is continuous with respect to $x\in{\mathbb R}$ because 
$M_{\lambda,\sigma_0}(t)$ is bounded. 
$\int_x^\infty e^{-(\sigma-\sigma_0)t}dF(t)$ 
is monotonic decreasing in $x$. That is, a decreasing function 
$\int_x^\infty e^{-(\sigma-\sigma_0)t}dF(t)$ 
is bounded above by a continuous function 
$e^{\sigma_0x}\int_0^\infty M_{\lambda,\sigma_0}(t-x)e^{-\sigma t}dF(t)$ 
outside a countable set of $x$, hence, so is for all $x\in{\mathbb R}$. 
Therefore, 
\begin{eqnarray}
e^{-\sigma_0 x}\int_x^\infty e^{-(\sigma-\sigma_0)t}dF(t)
\le\int_0^\infty M_{\lambda,\sigma_0}(t-x)e^{-\sigma t}dF(t)
\label{eqn:theo6:1}
\end{eqnarray}
holds for all $x\in{\mathbb R}$.

Since $M_{\lambda,\sigma_0}$ is of type $-2\sigma_0\lambda$, 
Paley-Wiener's theorem \cite{rud} shows that the support 
of its Fourier transform $\hat{M}_{\lambda,\sigma_0}$ 
is contained in the interval $[2\sigma_0\lambda,-2\sigma_0\lambda]$. 
Writing $\Lambda=-2\sigma_0\lambda$, we have by the 
inverse Fourier transform 
\begin{eqnarray}
M_{\lambda, \sigma_0}(t)=\frac{1}{2\pi}\int_{-\Lambda}^\Lambda
\hat{M}_{\lambda, \sigma_0}(\tau)e^{it\tau}d\tau,
\end{eqnarray}
and by substituting $t\to t-x$, $\tau\to-\tau$, 
\begin{eqnarray}
M_{\lambda, \sigma_0}(t-x)=\frac{1}{2\pi}\int_{-\Lambda}^\Lambda
\hat{M}_{\lambda, \sigma_0}(-\tau)e^{-i(t-x)\tau}d\tau.
\end{eqnarray}
Then we have by Fubini's theorem
\begin{align}
\frac{1}{x}\int_0^\infty M_{\lambda, \sigma_0}(t-x)
e^{-\sigma t}dF(t)
&=\frac{1}{2\pi x}\int_{-\Lambda}^\Lambda \hat{M}_{\lambda, \sigma_0}(-\tau)
e^{ix\tau}d\tau\int_0^\infty e^{-(\sigma+i\tau)t}dF(t)\\
&=\frac{1}{2\pi x}\int_{-\Lambda}^\Lambda \hat{M}_{\lambda, \sigma_0}(-\tau)
\varphi(\sigma+i\tau)e^{ix\tau}d\tau.\label{eqn:theo6:2}
\end{align}
Since $s=\sigma_0$ is supposed to be a pole of order 2, 
the principal part $\varphi_0(s)$ of $\varphi(s)$ at $s=\sigma_0$ 
is given as 
\begin{eqnarray}
\varphi_0(s)=\sum_{j=1}^2\frac{A_j}{(s-\sigma_0)^j},\ A_2\neq 0.
\end{eqnarray}
The inverse Laplace transform of $\varphi_0(s)$ is 
\begin{eqnarray}
\sum_{j=1}^2A_jt^{j-1}e^{\sigma_0 t},\ t\ge 0,
\end{eqnarray}
then, in a similar way as (\ref{eqn:theo6:2}), we have
\begin{align}
\frac{1}{x}\int_0^\infty M_{\lambda, \sigma_0}(t-x)e^{-\sigma t}
\sum_{j=1}^2A_jt^{j-1}e^{\sigma_0 t}dt
&=\frac{1}{2\pi x}\int_{-\Lambda}^\Lambda\hat{M}_{\lambda, \sigma_0}
(-\tau)\sum_{j=1}^2\frac{A_j}{(\sigma+i\tau-\sigma_0)^j}e^{ix\tau}d\tau.
\label{eqn:theo6:3}
\end{align}
We write $\xi(s)=\varphi(s)-\varphi_0(s)$, then 
$\xi(s)$ is analytic in a neighborhood of $s=\sigma_0$. 
By subtracting (\ref{eqn:theo6:3}) from (\ref{eqn:theo6:2}), 
\begin{align}
\frac{1}{x}\int_0^\infty M_{\lambda, \sigma_0}(t-x)
e^{-\sigma t}dF(t)
&=\sum_{j=1}^2\frac{A_j}{x}\int_0^\infty M_{\lambda, \sigma_0}(t-x)
t^{j-1}e^{-(\sigma-\sigma_0)t}dt\nonumber\\
&\hspace{5mm}+\frac{1}{2\pi x}\int_{-\Lambda}^\Lambda
\hat{M}_{\lambda, \sigma_0}(-\tau)\xi(\sigma+i\tau)e^{ix\tau}d\tau.
\label{eqn:theo6:4}
\end{align}
Then from (\ref{eqn:theo6:1}), (\ref{eqn:theo6:4}) we have
\begin{align}
\frac{1}{x}e^{-\sigma_0 x}\int_x^\infty e^{-(\sigma-\sigma_0)t}dF(t)
&\le\sum_{j=1}^2\frac{A_j}{x}\int_0^\infty M_{\lambda, \sigma_0}(t-x)
t^{j-1}e^{-(\sigma-\sigma_0)t}dt\nonumber\\
&\hspace{5mm}+\frac{1}{2\pi x}\int_{-\Lambda}^\Lambda
\hat{M}_{\lambda, \sigma_0}(-\tau)\xi(\sigma+i\tau)e^{ix\tau}d\tau.
\end{align}
By taking $\omega>0$ sufficiently large, or equivalently, 
taking $\lambda>0$ sufficiently small, 
$\xi(\sigma+i\tau)\to\xi(\sigma_0+i\tau)$ as 
$\sigma\downarrow\sigma_0$ uniformly in $\tau\in[-\Lambda,\Lambda]$ 
with $\Lambda=-2\sigma_0\lambda$. Hence by the dominating convergence 
theorem, we have 
\begin{align}
\frac{1}{x}e^{-\sigma_0 x}P(X>x)
&\le\sum_{j=1}^2\frac{A_j}{x}\int_0^\infty M_{\lambda, \sigma_0}(t-x)
t^{j-1}dt
+\frac{1}{2\pi x}\int_{-\Lambda}^\Lambda\hat{M}_{\lambda, \sigma_0}
(-\tau)\xi(\sigma_0+i\tau)e^{ix\tau}d\tau.\label{eqn:theo6:5}
\end{align}
When $x$ tends to infinity, the first term of the right 
hand side of (\ref{eqn:theo6:5}) becomes
\begin{align}
\lim_{x\to\infty}\sum_{j=1}^2\frac{A_j}{x}
\int_0^\infty M_{\lambda, \sigma_0}(t-x)t^{j-1}dt
=A_2\int_{-\infty}^\infty M_{\lambda, \sigma_0}(t)dt,
\end{align}
while the second term tends to $0$ by the Riemann-Lebesgue 
theorem. Hence from (\ref{eqn:theo6:5}) we have
\begin{eqnarray}
\limsup_{x\to\infty}\frac{1}{x}e^{-\sigma_0 x}P(X>x)
\le A_2\int_{-\infty}^\infty M_{\lambda, \sigma_0}(t)dt
<\infty.\nonumber
\end{eqnarray}

In a similar way, we have for minorant $m_\omega(t)$
\begin{eqnarray}
\liminf_{x\to\infty}\frac{1}{x}e^{-\sigma_0 x}P(X>x)
\ge A_2\int_{-\infty}^\infty m_{\lambda, \sigma_0}(t)dt
>0.\nonumber
\end{eqnarray}

\medskip

\noindent{\bf Proof of Theorem \ref{theo:3}:} From Theorem \ref{theo:6}, 
writing $C_1=A_2\int_{-\infty}^\infty m_{\lambda,\sigma_0}(t)dt>0$, \\
$C_2=A_2\int_{-\infty}^\infty M_{\lambda,\sigma_0}(t)dt>0$, we see 
that for any $\epsilon>0$ there exists $x_0$ with
\begin{align}
C_1-\epsilon<x^{-1}e^{-\sigma_0x}P(X>x)<C_2+\epsilon,\ 
\forall x\geq x_0.\label{eqn:proof3}
\end{align}
By taking logarithm of each side of (\ref{eqn:proof3}), 
we have
\begin{align}
\ds\lim_{x\to\infty}\frac{1}{x}\log P(X>x)=\sigma_0.
\end{align}

\bigskip

We immediately obtain the following 
corollary of Theorem \ref{theo:3}.

\medskip

\begin{coro} Let $X$ be a random variable taking 
non-negative integral values, and \\
$f(z)=\sum_{n=0}^\infty P(X=n)z^n$ the probability 
generating function of $X$. The radius of convergence 
of $f(z)$ is denoted by $r$ and $1<r<\infty$ is 
assumed. If $z=r$ is a pole of $f(z)$, then the 
tail probability $P(X>n)$ decays exponentially.
\end{coro}

\medskip

\section{Lemmas and Theorems for a Pole of Arbitrary Order}
In the last section, we proved lemmas and theorems for a 
pole of order 2. We here provide the statement of lemmas 
and theorems for a pole of arbitrary order, corresponding 
to each lemma and theorem in the last section. 

Let $D$ be the order of the pole $s=\sigma_0$ of the 
Laplace-Stieltjes transform $\varphi(s)$. First, define 
\begin{eqnarray}
K=\left\{\begin{array}{ll}
D, & \mbox{\rm if $D$ is odd},\\
D+1, & \mbox{\rm if $D$ is even},
\end{array}\right.
\end{eqnarray}
so, $K$ is always odd. Then, define 
\begin{align}
R(v)&=\frac{v^K}{1-e^{-v}},\ v\in\mathbb{R},\\[2mm]
r(v)&= v^K,\ v\in{\mathbb R},\\
\tilde{R}_\omega(v)&= R(v+\omega)
-\sum_{k=0}^{K-1}\frac{R^{(k)}(\omega)}{k!}v^k,\ \omega>0,\ v\in{\mathbb R},
\end{align}
where $R^{(k)}$ denotes the $k$-th derivative of $R$.

We show below Lemmas \ref{lem:1}$^\ast$-\ref{lem:5}$^\ast$ 
and Theorem \ref{theo:6}$^\ast$ corresponding to Lemmas 
\ref{lem:1}-\ref{lem:5} and Theorem \ref{theo:6} in 
the last section, respectively. 

\medskip

\noindent{\bf Lemma 1$^\ast$} {\it There exists $\omega_0>0$ such that 
for any $\omega\ge\omega_0$}, 
\begin{eqnarray}
\left\{\begin{array}{ll}
\tilde{R}_\omega(v)\ge 0, & v\ge 0,\\[2mm]
\tilde{R}_\omega(v)\le 0, & v<0.
\end{array}\right.
\end{eqnarray}

\medskip

\noindent{\bf Lemma 2$^\ast$} {\it For sufficiently large $v>0$},
\begin{eqnarray}
R^{(k)}(v)<\ds\frac{r^{(k)}(v)}{1-e^{-v}},\ k=1,3,\cdots,K.
\end{eqnarray}

\medskip

\noindent{\bf Lemma 3$^\ast$} {\it There exists $\omega_0>0$ such that 
for any $\omega\ge\omega_0$}, 
\begin{eqnarray}
\tilde{R}_\omega(v)\le\frac{1}{1-e^{-\omega}}v^K,\ v\ge 0.
\end{eqnarray}

\medskip

\noindent{\bf Lemma 4$^\ast$} {\it There exists $\omega_0>0$ such that 
for any $\omega\ge\omega_0$}, 
\begin{eqnarray}
-\int_{-\infty}^0\tilde{R}_\omega(v)e^{tv}dv
<\frac{R^{(K)}(\omega)}{t^{K+1}},\ t>0.
\end{eqnarray}

\bigskip

Define two functions $M_\omega(t)$ and $m_\omega(t)$ 
as follows.
\begin{align}
M_\omega(t)&=\frac{1}{K!}\left(\frac{\sin\pi t}{\pi}\right)^{K+1}
Q_\omega(t),\ -\infty<t<\infty,\\
Q_\omega(t)&=\sum_{n=0}^\infty e^{-n\omega}\sum_{k=1}^{K+1}
\frac{(-1)^kr^{(k-1)}(\omega)}{(t-n)^k}
+\sum_{k=1}^K\frac{(-1)^{k-1}R^{(k-1)}(\omega)}{t^k},\\[1mm]
m_\omega(t)&=M_\omega(t)-\ds\frac{1}{1-e^{-\omega}}
\left(\ds\frac{\sin\pi t}{\pi t}\right)^{K+1},\ -\infty<t<\infty.
\end{align}

\medskip

\noindent{\bf Lemma 5$^\ast$} {\it There exists $\omega_0>0$ such that 
for any $\omega\ge\omega_0$, 
$M_\omega(t)$ is a majorant for $E_\omega(t)$ of 
type $(K+1)\pi$ and $\int_{-\infty}^\infty M_\omega(t)dt<\infty$, 
moreover, $m_\omega(t)$ is a minorant for $E_\omega(t)$ 
of type $(K+1)\pi$ and $\int_{-\infty}^\infty m_\omega(t)dt>0$.}

\bigskip

Let $\lambda=2\pi/\omega$ and define 
\begin{align}
M_{\lambda, \sigma_0}(t)&=M_\omega(-\frac{\sigma_0}{\omega}t),
\ t\in{\mathbb R}, \\
m_{\lambda, \sigma_0}(t)&=m_\omega(-\frac{\sigma_0}{\omega}t),
\ t\in{\mathbb R}.
\end{align}
Both $M_{\lambda, \sigma_0}(t)$ and $m_{\lambda, \sigma_0}(t)$ 
are of type $-2\sigma_0\lambda$. By Lemma \ref{lem:5}$^\ast$, 
there exists $\lambda_0>0$ such that for any $\lambda<\lambda_0$, 
$M_{\lambda,\sigma_0}(t)$ is a majorant and $m_{\lambda,\sigma_0}(t)$ 
is a minorant for $E_{\sigma_0}(t)$, respectively. 

We can calculate 
$\int_{-\infty}^\infty M_{\lambda, \sigma_0}(t)dt$ and 
$\int_{-\infty}^\infty m_{\lambda, \sigma_0}(t)dt$ as 
follows.
\begin{align}
\int_{-\infty}^\infty M_{\lambda, \sigma_0}(t)dt
&=-\frac{\omega}{\pi\sigma_0(1-e^{-\omega})}
\int_{-\infty}^\infty \left(\frac{\sin t}{t}\right)^{K+1}dt\nonumber\\
&\hspace*{7mm}-\frac{\omega}{K!\sigma_0}\sum_{k=2,\, k\mbox{\rm\scriptsize:even}}^{K-1}
\frac{1}{\pi^{K+2-k}}\left\{\frac{r^{(k-1)}(\omega)}{1-e^{-\omega}}
-R^{(k-1)}(\omega)\right\}
\int_{-\infty}^\infty\frac{\sin^{K+1}t}{t^k}dt,\\
&<\infty,
\end{align}
\begin{align}
\int_{-\infty}^\infty m_{\lambda, \sigma_0}(t)dt
&=-\frac{\omega}{K!\sigma_0}\sum_{k=2,\, k\mbox{\rm\scriptsize:even}}^{K-1}
\frac{1}{\pi^{K+2-k}}\left\{\frac{r^{(k-1)}(\omega)}{1-e^{-\omega}}
-R^{(k-1)}(\omega)\right\}
\int_{-\infty}^\infty\frac{\sin^{K+1}t}{t^k}dt\\
&>0.
\end{align}

\noindent{\bf Theorem \ref{theo:6}$^\ast$} {\it Let $X$ be a non-negative 
random variable and $\varphi(s)$ be the Laplace-Stieltjes 
transform of the probability distribution function 
$F(x)$ of $X$. The abscissa of convergence of $\varphi(s)$ 
is denoted by $\sigma_0$ and $-\infty<\sigma_0<0$ is 
assumed. Let $s=\sigma_0$ be a pole of $\varphi(s)$ of 
order $D$ and $\varphi(s)$ be analytic on the interval 
$\{s=\sigma_0+i\tau\,|\,2\sigma_0\lambda<\tau<-2\sigma_0\lambda\}$ 
for some $\lambda>0$ except $s=\sigma_0$. Write $A_D$ as 
the coefficient of $(s-\sigma_0)^{-D}$ in the Laurent 
expansion of $\varphi(s)$ at $s=\sigma_0$. Then, we have}
\begin{align}
A_D\int_{-\infty}^\infty m_{\lambda, \sigma_0}(t)dt
&\le\liminf_{x\to\infty}x^{-D+1}e^{-\sigma_0x}P(X>x)\nonumber\\
&\le\limsup_{x\to\infty}x^{-D+1}e^{-\sigma_0x}P(X>x)\nonumber\\
&\le A_D\int_{-\infty}^\infty M_{\lambda, \sigma_0}(t)dt.\nonumber
\end{align}

\section{Conclusion}
We investigated a sufficient condition for the exponential 
decay of the tail probability $P(X>x)$ of a non-negative 
random variable $X$. For the Laplace-Stieltjes transform 
$\varphi(s)$ of the probability distribution function of 
$X$ with abscissa of convergence $\sigma_0,\ -\infty<\sigma_0<0$, 
we showed that if $s=\sigma_0$ is a pole of $\varphi(s)$ 
then the tail probability decays exponentially. If 
$\varphi(s)$ is given explicitly, then this sufficient 
condition is easy to check. For the proof of our 
main theorem, we extended Graham-Vaaler's Tauberian 
theorem to the case that the order of the pole of 
$\varphi(s)$ is arbitrary. Now, I have proved the conjecture 
that was written in \cite{nak2}.

\end{document}